\documentclass[reqno,10pt]{article}
\usepackage{hyperref}
\usepackage{amsmath}
\usepackage{amssymb}
\usepackage{mathrsfs}

\newtheorem{thm}{Theorem}[section]
\newtheorem{lemma}[thm]{Lemma}

\newtheorem{corol}[thm]{Corollary}
\newtheorem{propos}[thm]{Proposition}
\newtheorem{rema}{Remark}[section]
\def\bp{\begin{propos}}
\def\ep{\end{propos}}
\def\bt{\begin{thm}}
\def\et{\end{thm}}
\def\bco{\begin{corol}}
\def\eco{\end{corol}}
\def\bl{\begin{lemma}}
\def\el{\end{lemma}}
\def\br{\begin{rema}}
\def\er{\end{rema}}
\def\be{\begin{equation}}
\def\ee{\end{equation}}
\def\ba{\begin{array}}
\def\ea{\end{array}}
\def\bena{\begin{eqnarray}}
\def\eena{\end{eqnarray}}

\def\P{{\mathbb P}}

\def\R{{\mathbb R}}
\def\Z{{\mathbb Z}}

\def\V{{\rm Vol}}
\def\1{I}
\def\imath{\textbf{i}}
\def\jmath{\textbf{j}}

\def\fC{{\cal C}}
\def\fD{{\cal D}}

\def\fE{{e}}

\def\hB{{\mathscr B}}

\def\fL{{\cal L}}

\def\a{{\alpha}}

\def\b{\binom}

\def\QED{\hfill$\square$\vskip 3mm}

\def\Dp{\displaystyle}
\def\Df{\Dp\frac}

\def\({\left(}
\def\){\right)}
\def\ln{\lg}

\begin{document}
\title{On The Modified Newman-Watts Small World\\ and Its Random Walk
\\[5mm]
\footnotetext{AMS classification (2000): Primary 05C 80; secondary 60K 35, 60K 37 } 
\footnotetext{Research supported in part by the Natural Science Foundation of China (under
grants 11271356, 11471222) and the Foundation
of Beijing Education Bureau (under grant KM201510028002) }}

\author{Xian-Yuan Wu,\ \ Rui Zhu}\vskip 10mm
\date{}
\maketitle

{\vskip-20mm
\begin{center}\begin{minipage}{10cm}

{\small \noindent School of Mathematical Sciences, Capital Normal
University, Beijing, 100048, China.\ \ Email: \texttt{wuxy@mail.cnu.edu.cn};\ \ \texttt{lmozi999@163.com};\\[-5mm]}

\end{minipage}
\end{center}
\vskip 2mm
\begin{center}
\begin{minipage}{10cm}
{\bf Abstract}:
{\small It is well known that adding ``long edges (shortcuts)" to a regularly constructed graph will make the resulted model a small world. Recently, \cite{W} indicated that, among all long edges, those edges with length proportional to the diameter of the regularly constructed graph may play the key role. In this paper, we modify the original Newman-Watts small world by adding only long special edges to the $d$-dimensional lattice torus (with size $n^d$) according to \cite{W}, and show that the diameter of the modified model and the mixing time of random walk on it grow polynomially fast in $\ln n$.


}
\end{minipage}
\end{center}}

\vskip 5mm
\section{Introduction and statement of the results}
\renewcommand{\theequation}{1.\arabic{equation}}
\setcounter{equation}{0}
{\it Small world effect}, the fact that the diameters of most networks are considerably smaller than their sizes,
is one of the most important features of real-world complex networks. General speaking, for the graph $G=(V,E)$ as a network with $V$ large enough, say $G$ exhibits the small world effect, if the diameter of $G$ is at most polynomially large in $\ln |V|$. In 1929, the existence of small world effect had been speculated upon in a remarkable short story by Karinthy \cite{K}. In 1960s, Milgram \cite{M,TM} carried out his famous ``small-world" experiments, which finally led to the popular concept of the ``six degrees of separation" \cite{G}. Recent influential studies on small world effect perhaps started with the work of Watts and Strogatz published in 1998 \cite{WS}. From then on, people were much more interested in studying the structure features (including {\it small world effect}, {\it scale-free property} and {\it navigability}, etc.) of complex networks. Nowadays, the small world effect has been studied and verified directly in a large number of different networks, see \cite[Table 3.1]{N} and the references therein.

To seek the underlying causes which make most networks small worlds, many models have been introduced and studied by physicists and mathematicians, see \cite{BC,BR,NW,WS} {\it etc.}. Actually, many models were introduced to reveal such a fact that adding ``long edges" to a regularly constructed (lattice-like) graph will make the resulted graph a small world, and we will call it the {\it adding-long-edges} mechanism. Recently, Wu introduced a new model in \cite{W}, the {\it Poisson Geometry Small World} (PGSW) model, which is obtained by adding special long edges to the largest cluster of the supercritical Poisson continuum percolation on the $d$-dimensional torus. The results of \cite{W} indicated that, for the {\it adding-long-edges} mechanism, among all added long edges, those edges with length proportional to the diameter of the original model may play the key role. In other words, to make the resulted graph a small world, it seems that adding shorter edges is neither sufficient nor necessary.

The PGSW model was built on the largest cluster of the supercritical Poisson continuum percolation, and its structure is highly complicated. For the diameter of the model, the results obtained in \cite{W} are non-optimal, furthermore, for different parameters, \cite{W} can only obtained lower or supper bound respectively. In other words, \cite{W} failed to provide such a model with its diameter growing polynomially fast in $\ln n$.
In the present paper, we will introduce a relatively simple model, it is called the {\it modified Newman-Watts model}, and then study its diameter. We hope we can obtain more precise results, at least, we can give good lower and supper bounds to the diameter simultaneously. Thus, we get an example to show that, under the special ``adding-long-edges" mechanism proposed in \cite{W}, the diameter of the network with size $n^d$ may really grow polynomially fast in $\ln n$.

Note that the original NW small world was introduced in \cite{NW}, it was obtained from $T_n$, the 1-dimensional lattice ring with $n$ sites, by adding a Poisson number of shortcuts with mean $n\rho/2$, and attaching then to randomly chosen pairs of sites. The original NW small world was also studied by Durrett in \cite{D}, where Durrett gave his setting as the following: let $\xi_x$; $x\in T_n$ is a sequence of i.i.d. Poisson random variables with mean $p$, let $\xi_x$ is the number of added shortcuts concerning vertex $x$. For any added shortcut, the other endvertex is chosen from $T_n$ independently and uniformly at random.

Now, we begin to introduce our modified NW model. Let $T_n^d$, $d\geq 1$, denote the $d$-dimensional lattice torus obtained from $B(n)=\Z^d\cap [0,n]^d$ by cohering its opposite faces. Let $d^T_\infty(\cdot,\cdot)$ denote the $l_\infty$ metric on $T^d_n$ inherited from the usual $l_\infty$ metric $d_\infty(\cdot,\cdot)$ on $\R^d$ defined by $d_\infty(x,y):=\Dp\max_{1\leq s\leq d}|x_s-y_s|$ for any $x,y\in\R^d$. For any given constants $\a$, $\beta$, $\sigma$ and $\zeta$ satisfying $0 <\a <\beta < 1/2$, $\sigma>0$ and $\zeta \in \R$, we construct a random graph $G_n = G_n(\a,\beta, \sigma,\zeta)$ from $T_{n}^{d}$ as the following: for any $u, v\in T_{n}^{d}$, if $\a n\leq d^T_\infty(u,v)\leq \beta n$, then we connect $u$ and $v$ independently by a ``long edge" with
probability
 \be\label{1.1} p_n=\sigma n^{-d}\ln^{\zeta}n;\ee
otherwise, we do nothing.
Denote by $G_{n}=(V_n, E_n)$ the resulted graph after all random long edges are added
to $T^{d}_{n}$. We then call $G_{n}$ the Modified NW Model. Note that $V_n$ is same as the vertex set of lattice torus $T^d_n$, $E_n$ is a random edge set which contains the edge set of $T^d_n$ as its subset.

It is straightforward to see that there are many differences between the two models. First, and most importantly, in the modified NW model, only long edges with length
in order $n$ are added, this obeys the essence of the special {\it adding-long-edges} mechanism suggested by \cite{W}. Second, in the modified model, no loop and double edge is added, this seems better to fit the situations of the real world networks. Third, the modified model deals with the high dimensional cases ($d \geq 1$). Finally, the modified model possesses more geometry, for example, for some subset $S \subset T^{d}_{n}$, $d=1$, let ${\cal L}(S; S^c)$ denote the random number of long edges between $S$ and $S^c$. In the original NW model, it's a Poisson number with mean $p|S| - p|S|^2/n$ (according to Durrett's setting); but in the modified model, the distribution of ${\cal L}(S; S^c)$ depends not only on the size ($|S|$) but also on the shape of $S$.


This paper will first study the diameter of $G_n$. Recall that in a graph $G$, the distance $\fD_G(u,v)$ between two vertices $u$ and $v$ is the length (number of edges) of the shortest path between them, and the {\it diameter} diam$(G)$ of a connected graph $G$ is the maximum distance between two vertices.

\bt\label{th1}
A: Suppose $0<\a<\beta<1/2$ with $\Gamma:=(2\beta)^{d}-(2\a)^{d}>1/2$. Then
\begin{description}
\item[](i) for $0<\zeta<1$ and $\sigma>0$ or, for $\zeta=0$ and $\sigma>0$ large enough, there exists constant $C_1>0$ such that
 \be\label{1.2}\lim_{n\rightarrow\infty}\P({\rm diam}(G_n)\leq C_{1}\ln^{3} n)=1;\ee
\item[](ii) for $\zeta>1$ and $\sigma>0$ or, for $\zeta=1$ and $\sigma>0$ large enough, there exists constant $C_2>0$ such that
 \be\label{1.3}\lim_{n\rightarrow\infty}\P({\rm diam}(G_n)\leq C_{2}\ln^{2} n)=1.\ee
 \end{description}

B: Suppose $\sigma>0$, $\zeta<1$ and $0<\a<\beta<1/2$. Then for any $0<\nu<(1-\zeta)/d$, we have
  \be\label{1.4}\lim_{n\rightarrow\infty}\P({\rm diam}(G_n)\geq \ln^{\nu} n)=1;\ee
   furthermore, if $\sigma>0$ is small enough, then
\be\label{1.5}\lim_{n\rightarrow\infty}\P({\rm diam}(G_n)\geq \ln^{(1-\zeta)/d} n)=1.\ee

\et

\br For $0\leq\zeta<1$, Theorem~\ref{th1} shows that the diameter of $G_n$, the network model (with size $n^d$) constructed by the special ``adding-long-edges" mechanism proposed in \cite{W}, may really grow polynomially fast in $\ln n$. \er

As Newman noted in \cite{N}, the ultimate goal of the study of the structure of networks is to understand and explain the workings of systems built upon those networks. Clearly, random walks on networks are just the simplest (but important) workings of systems built upon networks. At the present paper, we will next study the {\it mixing time} of random walk on $G_n$. For basic concepts on mixing time and related problems on mathematics and statistical physics, one may refer to \cite{LPW} and the references therein. For mixing time of random walk on complex networks, one may refer to \cite{D,LPW}.

In a graph $G=(V,E)$, for any $u,v\in V$, let $d_G(u)$ be the degree of $u$ in $G$, and write $u\sim v$ if $u$ and $v$ are neighbors in $G$. Let $\Delta(G)$ denote the {\it maximum degree} of $G$, i.e. $\Delta(G):=\max\{d_G(u):u\in V\}$. For any $S\subset V$, define $\V(S):=\sum_{u\in S}d_{G}(u)$.

Let's consider our modified NW model $G_n=(V_n,E_n)$, where $V_n=T^d_n$. For any $u,v\in V_n$, we define a transition kernel by $P(u,u)= 1/2$, $P(u,v)=1/2d_{G_n}(u)$ if $u\sim v$
and $P(u,v)=0$ otherwise. A discrete time Markov chain $\{X_t:t\geq 0\}$ on $V_n$ with transition kernel $(P(u,v))$
is called the {\it lazy random walk} on $G_n$. Note that $\pi(u):=d_{G_n}(u)/D$ where $D=\sum_{v\in V_n}d_{G_n}(v)$, defines a reversible stationary distribution of $\{X_t\}$ since $\pi(u)P(u,v)=\pi(v)P(v,u)$. By the basic theory of Markov chains, for any initial state $u\in V_n$, the distribution of $X_t$, i.e. $P^t(u,\cdot):=\P(X_t\in\cdot\mid X_0=u)$, converges weakly to $\pi$ as $t\rightarrow\infty$. To measure convergence to equilibrium, we will use the
{\it total variation distance} $$||P^t(u,\cdot)-\pi||_{TV}:=\frac 12\sum_{v\in V_n}|P^t(u,v)-\pi(v)|.$$ The mixing
time of $\{X_t:t\geq 0\}$ is defined by \be\label{1.6} T_{\rm mix}:=\min\left\{t:\max_{u\in V_n}||P^t(u,\cdot)-\pi||_{TV}<1/e\right\}.\ee

The second result of the present paper is about $T_{\rm mix}$ and we state it as follows

\bt\label{th2}
A: Suppose $0<\a<\beta<1/2$ with $\Gamma:=(2\beta)^{d}-(2\a)^{d}>1/2$. Then
\begin{description}
\item[](i) for $0<\zeta<1$ and $\sigma>0$ or, for $\zeta=0$ and $\sigma>0$ large enough, there exists constant $C_3>0$ such that
 \be\label{1.7}\lim_{n\rightarrow\infty}\P(T_{\rm mix}\leq C_{3}\ln^{5} n)=1;\ee
\item[](ii) for $\zeta>1$ and $\sigma>0$ or, for $\zeta=1$ and $\sigma>0$ large enough, there exists constant $C_4>0$ such that
 \be\label{1.8}\lim_{n\rightarrow\infty}\P(T_{\rm mix}\leq C_{4}\ln n)=1.\ee
 \end{description}

B: Suppose $\sigma>0$, $\zeta<1$ and $0<\a<\beta<1/2$. Then for any $0<\nu<2(1-\zeta)/d$, we have
  \be\label{1.9}\lim_{n\rightarrow\infty}\P(T_{\rm mix}\geq \ln^{\nu} n)=1.\ee
\et

The rest of the paper is arranged as follows. In Section 2, we bound the maximum degree of $G_n$ from above, and bound the {\it isoperimetric constant} of $G_n$ and the {\it conductance} of random walk $\{X_t\}$ from below. Finally, we prove Theorems~\ref{th1},~\ref{th2} in Section 3.

\vskip 5mm
\section{Isoperimetric constant and conductance }
\renewcommand{\theequation}{2.\arabic{equation}}
\setcounter{equation}{0}

On the graph $G_n=(V_n,E_n)$, for any $v\in V_n$, let $\Lambda_n(u)=\{v\in V_n:\a n\leq d^T_\infty(u,v)\leq \beta n\}$.  For any vertex sets $S,\ S'\subset V_n$, let $\Lambda(S,S')$ denote the set of {\it unordered} vertex pairs $\{u,v\}$ with $u\in S$, $v\in S'$ and $v\in \Lambda_n(u)$. Let $N(S,S')=|\Lambda(S,S')|$. To any {\it unordered} vertex pair $\{u,v\}\in \Lambda(V_n,V_n)$, independently, we assign a random variable $J_{u,v}(=J_{v,u})$ satisfying $\P(J_{u,v}=1)=1-\P(J_{u,v}=0)=p_n$. For any vertex sets $S,S'\subset V_n$, define
$$\fL(S,S'):=\Dp\sum_{\{u,v\}\in \Lambda(S,S')}J_{u,v}.$$
Let $\fL(S):=\fL(S,S)+\fL(S,S^c)$. Clearly, $\fL(S,S)$ and $\fL(S,S^c)$ are independent binomial random variables with parameters $(N(S,S),p_n)$ and
$(N(S,S^c),p_n)$ respectively.
Let $N(S)=N(S,S)+N(S,S^c)$, then $\fL(S)$ is the binomial random variable with parameter $(N(S),p_n)$.

New, we introduce the following large deviation inequality for binomial distribution for future use.

\bl\label{le2.1}{Suppose $Z\sim b(n,p)$, the binomial random variable with parameter $(n,\ p)$. Then
\be\label{2.1}\P(Z\geq zn)\leq \exp(-I(z)n),\ z>p \ {\rm and }\ \ \P(Z\leq zn)\leq \exp(-I(z)n), \ z<p,\ee where $I(z)$ is the common rate function defined by \be\label{2.2} I(z):=z\ln\frac{zq}{(1-z)p}-\ln \frac q{1-z}, \ p\not=z\in (0,1),\ \ p+q=1.\ee} Especially for small $p$, (\ref{2.1}) can be rewritten as
\be\label{2.3}\ba{rl}&\P(Z\geq zpn)\leq \exp(-\gamma(z)pn)\ {\it for}\ z>1, \ {\rm and }\\[3mm]&\P(Z\leq zpn)\leq \exp\(-\Df 12\gamma(z)pn\)\ {\it for}\ 0<z<1,\ea\ee with $\gamma(z)=z\ln z-z+1$. When $z>1$ is large enough, the first inequality in (\ref{2.3}) can be rewritten as
\be\label{2.3'}\P(Z\geq zpn)\leq \exp(-zpn).\ee\el

{\it Proof.} (\ref{2.1}) follows from the proof of the classical Cram\'er's Theorem \cite{C}. (\ref{2.3}) follows from (\ref{2.1}) by using the Taylor's expansion of $I(zp)$ for small $p$.\QED

First of all, we shall bound  $\Delta(G_n)$, the {\it maximum degree} of $G_n$ from above. Actually, we have the following proposition.



\bp\label{le2.2} Suppose $\zeta\geq 0$. Then for some constant $D>0$,
\be\label{2.4}\lim_{n\rightarrow\infty}\P(\Delta(G_n)\leq D\ln^{\zeta \vee 1} n)=1.\ee\ep
{\it Proof.} For any $u\in V_n$, $d_{G_n}(u)=2d+b(\Gamma n^d,p_{n})$, where $\Gamma=(2\beta)^d-(2\a)^d$. Taking $D$ large enough and using the large deviation inequality (\ref{2.3}), we have $$\P(b(\Gamma n^d,p_{n})\geq D\ln^{\gamma}n)\leq n^{-\gamma(z)\Gamma\sigma\ln^{\zeta-1}n},$$ where $z={D\ln^{\gamma-\zeta}n}/{\sigma\Gamma}>1$ and $\gamma=\zeta \vee 1$. For large enough $D$ and $n$, we always have $\gamma(z)\Gamma\sigma\ln^{\zeta-1}n>d$, then the proposition follows immediately from the definition of $\Delta(G_n)$.\QED
%

For the lazy random walk $\{X_t:t\geq 0\}$ on $G_n$, let $Q(u,v):=\pi(u)P(u,v)$ and $Q(S,S^c):=\sum_{u\in S}\sum_{v\in S^c}Q(u,v)$. Define
$$h:=\min_{S:\pi(S)\leq 1/2}\frac{Q(S,S^c)}{\pi(S)}$$ to be the {\it conductance} of $\{X_t:t\geq 0\}$.
Letting $\fE(S,S^c)$ be the number of edges between $S$ and $S^c$, i.e.
$\fE(S,S^c)=|\{(u,v)\in{E_{n}}:u\in{S},v\in{S^c}\}|$, we have \be\label{2.5}h=\frac 12\min_{S:\pi(S)\leq 1/2}\frac{\fE(S,S^c)}{\V(S)},\ee where $\pi(S)=\sum_{u\in S}\pi(u)$.

Another interesting quality on $G_n$ is the {\it edge isoperimetric constant} $\iota$ defined by
\be\label{2.6}\iota:=\min_{S:|S|\leq |V_n|/2}\frac{\fE(S,S^c)}{|S|}.\ee
In the rest of this section, we will try to give lower bounds to $h$ and $\iota$. Using these lower bounds, we then finish the proofs of our main results in the next section.

In the following proofs, we will use the notations $\Omega(b_n)$, $O(b_n)$ and $\Theta(b_n)$. In fact, we use $a_n=\Omega(b_n)$, $a_n=O(b_n)$ and $a_n=\Theta(b_n)$ to denote $a_n\geq cb_n$, $a_n\leq Cb_n$ and $cb_n\leq a_n\leq C b_n$ respectively for some $0<c,\ C<\infty$.

\bl\label{le2.3}Suppose  $\zeta\geq 0$. Then for small enough $a>0$,
\be\label{2.7}\lim_{n\rightarrow\infty}\P\(\bigcap_{S:|S|\geq (1-a)n^{d}}\left\{\pi(S)>\frac 12\right\}\)=1.\ee \el
{\it Proof.} For any $S\subset V_n$, we have
\be\label{2.8}\pi(S)=\sum _{u\in S}\pi(u)\geq\Df{d|S|+2\fL(S,S)+\fL(S,S^c)}{dn^d+2\fL(V_n)},\ee
where $\fL(S)=\fL(S,S)+\fL(S,S^c)$.

On the one hand, for $S$ with $|S|\geq (1-a)n^{d}$ and $a>0$ small enough,
\be\label{2.9}\ba{rl}&\Df{d|S|}{dn^d}\geq\Df{2d(1-a)n^{d}}{2dn^{d}}=(1-a)>\Df{1}{2}.\\[2mm]\ea\ee
On the other hand, for any given $S$ with $|S|\geq(1-a)n^d$, by definitions, $\fL(S,S)\sim b(N(S,S),p_n)$, $\fL(S^c)\sim b(N(S^c), p_n)$ are independent binomial random variables, and for small enough $a>0$,
\be\label{2.10}N(S^c)\leq(1-a)n^{d}\cdot{an^{d}}+\Df{(an^d)^2}{2}=\((1-a)a+\Df{a^2}{2}\)\cdot{n^{2d}}\leq\frac A 5 n^{2d},\ee
\be\label{2.11}N(S,S)\geq\(\frac{\Gamma^2}2-(1-a)a-\Df{a^{2}}{2}\)n^{2d}\geq\frac {4A}5 n^{2d},\ee
with $A=\Gamma^2/2$. Let $Y_{1}\sim b(\frac A 5 n^{2d},p_n)$, $Y_{2}\sim b(\frac {4A} 5n^{2d},p_n)$ are two independent random variables, then $$\P\(\fL(S,S)<\fL(S^c)\)\leq\P(Y_{2}<Y_{1}).$$ Using the inequality (\ref{2.3}), we know that both $\P(Y_{2}<{\frac {3A} 4 n^{2d}p_{n}})$ and $\P(Y_{1}>\frac A 4{n^{2d}p_{n}})$ are less than $\exp\(-\Omega(n^d\ln^{\zeta} n)\)$, hence
\be\label{2.12'}\P\(\fL(S,S)<\fL(S^c)\)\leq \exp\(-\Omega(n^d\ln^{\zeta} n)\).\ee

Note that, by (\ref{2.9}) and the the fact that $\fL(S,S)\geq\fL(S^c)$,
\be\label{2.12}\ba{rl}\pi(S)&\geq\Df{d|S|+2\fL(S,S)+\fL(S,S^c)}{dn^d+2\fL(V_n)}\geq\Df{d|S|+2\fL(S,S)}{dn^d+2\fL(S,S)+2\fL(S^{c})}\\[4mm]
&>\Df{\frac 12 dn^d+\fL(S,S)+\fL(S^{c})}{dn^d+2\fL(S,S)+2\fL(S^{c})}=\Df{1}{2}.
\ea\ee
Then, by (\ref{2.12'}), we obtain
\be\label{2.13}\P\(\pi(S)>\frac 12\)\geq 1-\exp\(-\Omega(n^{d}\ln^{\zeta}n)\).\ee
Now, let $M_a:=|\{S\subset V_n:|S|\geq (1-a)n^d\}|$. To finish the proof of the lemma, it remains to bound $M_a$ from above.
By Lemma 6.3.3 in \cite{D}, the number of $S\subset V_n$ with $|S|=an^{d}$ is
\be\label{2.14}\ba{rl}\Dp\b{n^d}{an^d}\leq \(\Df{n^d\cdot e}{an^d}\)^{an^d}&\hskip-3mm=\exp\left\{an^d\left[\ln\(\frac{n^d}{an^d}\)+1\right]\right\}\\[3mm]&\hskip-3mm=\exp\{a(\ln(1/a)+1)n^{d}\}.\ea\ee
Then for small enough $a>0$, \be\label{2.15}M_a=\Dp\sum_{s=(1-a)n^d}^{n^d}\b{n^d}{s}=\Dp\sum_{s=1}^{an^d}\b{n^d}{s}\leq an^d\b{n^d}{an^d}\leq\exp\{O\(c(a)n^d\)\},\ee where $c(a)>0$ and tends to $0$ as $a\rightarrow 0$.
Combining (\ref{2.13}) and (\ref{2.15}), we obtain
$$\P\(\bigcup_{S:|S|\geq (1-a)|V_n|}\left\{\pi(S)\leq\Df {1}{2}\right\}\)\leq M_a\cdot\exp\(-\Omega(n^{d}\ln^{\zeta}n)\)\rightarrow 0\ \ {\rm as}\ \ n\rightarrow\infty,$$ for small enough $a>0$.\QED


Let \be\label{2.21}\hat \iota=\min_{S:|S|\leq(1-a)|V_n|}\frac{\fE(S,S^c)}{\V(S)}.\ee Then, by Lemma~\ref{le2.3},
\be\label{2.22}\P\(h\geq \hat \iota/2\)\rightarrow 1,\ \ {\rm as}\ n\rightarrow \infty.\ee So, to bound $h$ from below, it suffices to bound $\hat \iota$ from bellow. In fact, we have

\bp\label{p2.4}Suppose $0<\a<\beta<1/2$ with $\Gamma:=(2\beta)^d-(2\a)^d>1/2$. Then
 \begin{description}
 \item[] (i) for $0<\zeta<1$ and $\sigma>0$ or, for $\zeta=0$ and $\sigma>0$ large enough, there exists $C_{6}>0$ such that
          \be\lim_{n\rightarrow\infty}\P\(\hat \iota\geq C_{6}\ln^{-2} n\)=1;\ee
 \item[] (ii) for $\zeta>1$ and $\sigma>0$ or, for $\zeta=1$ and $\sigma>0$ large enough, there exists $C_{7}>0$ such that
          \be\lim_{n\rightarrow\infty}\P\(\hat \iota\geq C_{7}\)=1.\ee

\end{description}
\ep

{\it Proof.} This is the main part of our proofs. In fact, we will develop a more complicated version of the approach used in \cite[Theorem 6.6.1]{D} and \cite{W}.

Let
$$\ba{rl}&\hB_1:=\{S\subset V_n: 1\leq|S|\leq M\ln n\},\ {\rm and}\\[3mm]
         &\hB_2:=\{S\subset V_n: M\ln n<|S|\leq an^{d}\},\ {\rm and}\\[3mm]
         &\hB_3:=\{S\subset V_n:an^{d}<|S|\leq (1-a)n^{d}\},\ea$$
         where $M>0$ is a large constant and $a>0$ is given in Lemma~\ref{le2.3}. Item (i) of this proposition follows from the following Lemmas~\ref{le2.6}, \ref{le2.7} and ~\ref{le2.9}, item (ii) of the proposition follows from the following Lemmas~\ref{le2.8} and~\ref{le2.9}.\QED

\bl\label{le2.6} Suppose $0<\a<\beta<1/2$. For $0\leq\zeta<1$ and $\sigma>0$, there exists $C_{8}>0$ such that
\be\label{2.25}\lim_{n\rightarrow\infty}\P\left\{\Df{\fE (S,S^c)}{\V(S)}\geq\Df{C_{8}}{\ln^{2}n},\ \forall\ S\in\hB_1\right\}=1.\ee\el

{\it Proof.} For any $S\subset V_n$, we have $\fE(S,S^c)\geq 1$, $\V(S)=\sum_{u\in S}d_{G_n}(u)\leq2dM\ln n+2\fL(S)$.
Let's consider the random variable $\fL(S)\sim b(N(S),p_n)$. Clearly,
\be\label{2.26}N(S)\leq|S|\Gamma n^d=:N_1(S), \ N(S)\geq\frac 12|S|\Gamma n^d=:N_2(S),\ee
note that this indicates that $N(S)=\Theta(n^d|S|)$.

 For any $S\in\hB_1$ and any constant $D_1>0$, we have
$$\P\(\fL(S)> D_1\ln^{2} n\)\leq\P\(b(N_1(S),p_n)> D_1\ln^{2} n\).$$
By the large deviation inequality (\ref{2.3'}),
\be\label{2.27}\P\(b(N_1(S),p_n)> D_1\ln^{2} n\)\leq \exp\(-D_1\ln^{2}n \).\ee

By Lemma 6.3.3 in \cite{D}, the number of $S\subset V_n$ with $|S|=s$ is
\be\label{3.25'}\b{n^d}{s}\leq \(\Df{n^d\cdot e}{s}\)^s=\exp\left\{s\left[\ln\(\frac{n^d}s\)+1\right]\right\}=\exp\{O(s\ln n )\}.\ee
Then \be\label{2.28}|\hB_1|\leq\sum_{s=1}^{M\ln n}\b{n^d}{s}\leq M\ln n\b{n^d}{M\ln n}\leq \exp\(O\(\ln^{2} n\)\).\ee
Note that $\fL(S)\leq D_1\ln^{2} n$ implies
$$\Df{\fE(S,S^c)}{\V(S)}\geq\Df{1}{2dM\ln n+2\fL(S)}\geq\Df{C_{8}}{\ln^{2} n}.$$
Taking $D_1$ large enough, we finish the proof of the lemma by (\ref{2.27}) and (\ref{2.28}).\QED

\br For mixing time of random walk on the original NW model, Durrett \cite{D} obtained a supper bound $O(\ln^3 n)$. In fact, to bound the corresponding conductance, Durrett declared that the quantity $\inf_{S\in\hB_1}\frac {{\fE}(S,S^c)}{\V(S)}$ is bounded from below by $\Omega (\ln^{-1}n)$ (see \cite[page 174, line 10-11]{D}), while we obtained the bound $\Omega(\ln^{-2} n)$ in Lemma~\ref{le2.6}. It is regretful that we can not give a proof to his declaration at the present time. \er

\bl\label{le2.7} Suppose $0<\a<\beta<1/2$. For $0<\zeta<1$ and $\sigma>0$, or for $\zeta=0$ and $\sigma>0$ large enough, there exists $C_{9}>0$ such that
\be\label{2.29}\lim_{n\rightarrow\infty}\P\left\{\Df{\fE (S,S^c)}{\V(S)}\geq\Df{C_{9}}{\ln^{2} n},\ \forall\ S\in\hB_2\right\}=1.\ee\el
{\it Proof.} For any vertex subset $S$, let $G_{S}$ be the subgraph of the torus $T^d_n$ with vertex set $S$ and edge set consists of all edges which connect vertex pairs in $S$. Denote by $\fC(S)$ the number of connected components of $G_{S}$. Let $\hB^{\geq}_2:=\{S\in\hB_2:\fC(S)\geq{|S|}/{\ln n}\},$
                $\hB^{<}_2=\hB_2\setminus \hB^{\geq}_2.$

At first considering,  for any $ S\subset\hB^{\geq}_2$ and $D_2>0$ large enough, by the large deviation inequality (\ref{2.3'}), we have
$$\P(\fL(S)\geq D_2|S|\ln n)\leq\P\(b(N_1(S),p_n)\geq D_2|S|\ln n\)\leq \exp(-D_2|S|\ln n).$$
Then, by (\ref{3.25'}), we have
\be\label{2.30}\ba{rl}\P\(\Dp\bigcup_{S\in \hB^{>}_2}\{\fL(S)\leq D_2|S|\ln n\}^c\)&\leq\Dp\sum^{an^d}_{s=M\ln n}\b{n^d}{s}\exp\(-\(D_2s\ln n\)\)\\[3mm]&\leq an^d\exp\{-[\(D_2\ln^2 n\)-O\(\ln^2 n\)]\}\rightarrow 0,\ea\ee as $n\rightarrow\infty$.
Note that $\fL(S)\leq D_2|S|\ln n$ implies,
\be\label{2.28'}\Df{\fE (S,S^c)}{\V(S)}\geq\Df{|S|/\ln n}{2d|S|+2D_2|S|\ln n}\geq\Df{C_9}{\ln^{2}n}.\ee

For any $ S\in\hB^{<}_2$, let's consider the random variables $\fL(S,S^c)$, $\fL(S,S)$ and $\fL(S)$. Recall that $\fL(S,S^c)\sim b(N(S,S^c),p_n)$, $\fL(S,S)\sim b(N(S,S),p_n)$ and $\fL(S)\sim b(N(S),p_n)$.
By definition
$$N(S,S^c) \geq \Dp\sum_{u\in S}\(\Gamma n^d-|S|\)\geq \sum_{u\in S}\(\Gamma n^d-an^d\)\geq (\Gamma-a)n^d\cdot|S|\geq(1-\epsilon) N(S), $$
 where $\epsilon=a/\Gamma$. Then $$ N(S,S)=N(S)-N(S,S^c)\leq N(S)-(1-\epsilon) N(S)=\epsilon N(S). $$
Now, let $W_{1}\sim b((1-\epsilon) N(S),p_n)$, $W_{2}\sim b(\epsilon N(S),p_n)$ and suppose that $W_{1}$ and $W_{2}$ are independent,
then,
$$\ba{rl}&\P\(\fL(S,S^c)\geq \Df{\fL(S)}2\)=\P(\fL(S,S^c)\geq \fL(S,S))\\&\geq \P(W_1\geq W_2)\geq \P\(W_1\geq \Df {N(S)p_n}2\geq  W_2\).\ea$$
By inequality (\ref{2.3}), both $\P\(W_1\leq \frac{N(S)p_n}2\)$ and $\P\( W_2\geq \frac{N(S)p_n}2\)$ are less than $\exp(-\Omega (\sigma |S|\ln^{\zeta}n))$,
so \be\label{2.31}\P \(\fL(S,S^c)\geq \Df{\fL(S)}2\)\geq1-\exp(-\Omega (\sigma |S|\ln^{\zeta}n)).\ee
Using the inequality (\ref{2.3}) again, we obtain
\be\label{2.32}\ba{rl}\P\(\fL(S)\leq |S| \)\leq\P\(b\(N_2(S),p_n\)\leq |S|\)
\leq\exp(-\Omega (\sigma |S|\ln^{\zeta}n)).\ea\ee
Note that $\fL(S,S^c)\geq {\fL(S)}/2\ {\rm and\ }\fL(S)\geq |S|$ imply
\be\label{2.33}\Df{e(S,S^c)}{\V(S)}\geq \Df{\fL(S,S^c)}{2d|S|+2\fL(S)}
\geq\Df{\fL(S)/2}{2d\fL(S)+2\fL(S)}\geq C_9.\ee

For any $M\ln n\leq s\leq an^d$ and any $1\leq j\leq s/\ln n=:j_s$, let $\hB^{<}_{2,s}=\{S\in\hB^{<}_{2}:|S|=s\}$, $\hB^{<}_{2,s,j}=\{S\in\hB^{<}_{2,s}:\fC(S)=j\}$.
Now, we have to bound $|\hB^{<}_{2,s,j}|$ from above.

In the case of $d=1$, same as Durrett discussed in \cite[page 172, case 1b]{D},
\be\label{2.34}\ba{rl}|\hB^{\leq}_{2,s,j}|&\leq \Dp\b{n}{j}\b{s-1}{j-1}\leq\b{n}{j_{s}}\b{s}{j_{s}}\leq\(\Df{n{e}}{j_{s}}\)^{j_{s}}\(\Df{s{e}}{j_{s}}\)^{j_s}\\[4mm]
&=\exp\left\{ j_{s} \(\ln\Df{n}{j_{s}}+\ln\Df{s}{j_{s}}+2\)\right\}\\[4mm]
&\leq\exp\left\{\Df{s}{\ln n}\Dp\(\ln n+\ln\ln n+2-\ln M \)\right\}\leq e^{2s}.\ea\ee
In case of $d\geq 2$, the situation is more complicated. To get a supper bound for $|\hB^{<}_{2,s,j}|$, one may first choose a subset, say, $\{v_1,v_2,\ldots,v_j\}\subset V_n$. Then, for any $v_i$, look it as a seed and let it grow into a connected subgraph $G_{v_i}=(V_{v_i}, E_{v_i})$ of $T^d_n$ such that $G_{v_i}, i=1,\ldots,j$ forms the $j$ connected components of $G_S$ with $S=\cup_{i=1}^j V_{v_i}$. Let $k_{v_i}(x_i)$ be the number of connected subgraph of $T_n^d$ that contains the point of $v_i$ and has size $x_i$, then, by the following Lemma~\ref{le2.4}, we have
\be\label{2.35}\ba{rl}|\hB^{<}_{2,s,j}|&\leq\Dp\sum_{\{v_1,\cdots,v_j\}\subset V_n}\sum_{(x_1,\ldots,x_j):\sum x_i=s,x_i\geq 1}\prod_{i=1}^j k_{v_i}(x_i)\\[6mm]
&\leq\Dp\b{n^d}{j}\b{s-1}{j-1}e^{2dcs}\leq\b{n^d}{j_{s}}\b{s}{j_{s}}e^{2dcs}\leq e^{O(s)}.
\ea\ee
Using (\ref{2.31}), (\ref{2.32}), (\ref{2.34}) and (\ref{2.35}), we obtain
 \be\label{2.36}\ba{rl}&\P\(\Dp\bigcup_{S\in \hB^{<}_2}\{\fL(S,S^c)\geq {\fL(S)}/2\ {\rm and\ }\fL(S)\geq |S|\}^c\)\\[5mm]
 &\leq\Dp\sum_{s=M\ln n}^{an^d}j_s\exp{(O(s))}\exp(-\Omega (\sigma s\ln^{\zeta}n))\\[5mm]
 &\leq\Dp\sum_{s=M\ln n}^{an^d}\Df{s}{\ln n}\exp\left\{-\(\Omega(\sigma s\ln^{\zeta}n)-O(s)\)\right\}\rightarrow 0,\ {\rm as} \ n\rightarrow\infty.\ea\ee
 Note that in the last line of (\ref{2.36}), when $\zeta>0$, $\sigma$ can be taken as any constant; but in case of $\zeta=0$, it requires $\sigma$ large enough.

The lemma follows immediately from (\ref{2.30}), (\ref{2.28'}), (\ref{2.33}) and (\ref{2.36}).\QED

\bl\label{le2.4} For any connected graph $G=(V,E)$, let $\Delta(G)$ denote its maximum degree and $k(x)$ be the number of connected subgraph that contains the point of $o\in V$ with size $x$. Then there exists constant $c>0$, such that
$$k(x)\leq e^{c\Delta(G)x}.$$\el
{\it Proof.} See the arguments in \cite[Eq: (4.24), page 81]{G}.\QED

\bl\label{le2.8} Suppose $0<\a<\beta<1/2$. For $\zeta>1$ and $\sigma>0$, or for $\zeta=1$ and $\sigma>0$ large enough, there exists  $C_{10}>0$ such that
\be\label{2.37}\lim_{n\rightarrow\infty}\P\left\{\Df{\fE
(S,S^c)}{\V(S)}\geq C_{10},\ \forall\ S\in\hB_1
\cup\hB_2\right\}=1.\ee\el
{\it Proof.} First, we point out that, for any $S\subset V_n$, (\ref{2.31}) and (\ref{2.32}) hold for
random variables $\fL(S; S^c)$ and $\fL(S)$.

Second, by using (\ref{3.25'}), (\ref{2.31}) and (\ref{2.32}), we obtain
 \be\label{2.42}\ba{rl}&\P\(\Dp\bigcup_{S\in\hB_1 \cup\hB_2}\{\fL(S,S^c)\geq {\fL(S)}/2\ {\rm and\ }\fL(S)\geq |S|\}^c\)\\[5mm]
 &\leq\Dp\sum_{s=1}^{an^d}\b{n^d}{s}\exp(-\Omega (\sigma s\ln^{\zeta}n))\\[5mm]
 &\leq a n^d\exp\(-\(\Omega\(\sigma\ln^{\zeta}n\)-O\(\ln n\)\)\)\rightarrow 0,\ {\rm as} \ n\rightarrow\infty.\ea\ee
Then, the lemma follows immediately from (\ref{2.33}) and (\ref{2.42}).\QED

\bl\label{le2.9} Suppose $0<\a<\beta<1/2$ with $\Gamma:=(2\beta)^d-(2\a)^d>1/2$. Then, for $\zeta>0$ and $\sigma>0$, or for $\zeta=0$ and $\sigma>0$ large enough, there exists $C_{11}>0$ such that
\be\label{2.43}\lim_{n\rightarrow\infty}\P\left\{\Df{\fE (S,S^c)}{\V(S)}\geq C_{11},\ \forall\ S\in\hB_3\right\}=1.\ee\el

{\it Proof.} The proof of this lemma is very similar to the proof of Lemmas~\ref{le2.7} and \ref{le2.8}. But just in this step, we have to use the condition of $\a$ and $\beta$: $(2\beta)^d-(2\a)^d>1/2.$

Suppose that $a$ is taken small enough such that
$$\Gamma-\frac 12\geq a.$$

For any $S\in \hB_3$, if $an^d\leq |S|\leq {n^d}/{2}$, then
$$\ba{rl}N(S,S^c)&=\Dp\sum_{u\in S}\sum_{v\in S^c\cap\Lambda_n(u)}1\geq\Dp\sum_{u\in S}\(\Gamma n^d-|S|\)=|S|\(\Gamma n^d-|S|\)\\[4mm]&\geq an^d\(\Gamma n^d-a n^d\)\geq a(\frac12-a)n^{2d}.\ea$$
Note that, in the first inequality we have used the condition $\Gamma>1/2$ to guarantee that all terms in the summation in the right hand side are nonnegative. The second inequality comes from the fact that the function $g(x)=x(1-x)$ in interval $[a/\Gamma,1/2\Gamma] $ takes its minimum at $x=a/\Gamma$.

On the other hand, if $an^d\leq |S^c|\leq {n^d}/{2}$, then
$$\ba{rl}N(S,S^c)&=N(S^c,S)=\Dp\sum_{u\in S^c}\sum_{v\in S\cap\Lambda_n(u)}1\geq\Dp\sum_{u\in S^c}(\Gamma n^d-|S^c|) \\[3mm]
&\geq an^d\(\Gamma n^d-an^d\)\geq a(\frac12-a)n^{2d}.\ea$$
So, for any $S\in\hB_3$, $$N_1(S)=\Gamma n^d|S|\leq (1-a)\Gamma n^{2d}\leq\frac{(1-a)\Gamma}{a(\frac12-a)}N(S,S^c).$$
Let $f(a)=a(\frac12-a)/(1-a)\Gamma$, then$$N(S,S^c)\geq f(a)N_1(S)\geq f(a)N(S).$$
Now, by the large deviation inequality (\ref{2.3}), we have
$$\ba{rl}&\P\(\fL(S,S^c)\leq \Df{1}{2}f(a)N(S)p_n\)\leq\Dp\P\(\fL(S,S^c)\leq \frac 12N(S,S^c)p_n\)\\[4mm]&\Dp\leq \exp\(-\frac12\gamma\(\frac 12\)N(S,S^c)p_n\)\leq\exp\(-\frac12\gamma\(\frac 12\)f(a)N(S)p_n\)\\[4mm]&\leq\exp(-\Omega(\sigma |S|\ln^{\zeta} n))\ea$$ and
$$\ba{rl}\P\(\fL(S)\geq 2N(S)p_n\)&\leq \exp\(-\gamma(2)N(S)p_n\)\leq  \exp\(-\gamma(2)N_2(S)p_n\)\\[2mm]
&\leq\exp(-\Omega(\sigma|S|\ln^{\zeta} n)).\ea$$

If $\fL(S,S^c)\leq f(a)N(S)p_n/2$ and $\fL(S)\geq 2N(S)p_n$, then
$\fL(S,S^c)\geq f(a)\fL(S)/4$. So
\be\label{2.48}\ba{rl}&\P\(\fL(S,S^c)\geq \Df {1}{4} f(a)\fL(S)\)\\[3mm]&\geq\P\(\fL(S,S^c)\geq \frac12f(a)N(S)p_n\ {\rm and\ }\fL(S)\leq 2N(S)p_n\)\\[3mm]&\geq 1-\exp\(-\Omega\(\sigma|S|\ln^{\zeta}n\)\).\ea\ee
Similar to (\ref{2.33}),  $\fL(S,S^c)\geq f(a)\fL(S)/4\ $ and $\fL(S)\geq |S|$ imply
\be\label{2.48'}\Df{e(S,S^c)}{\V(S)}\geq \Df{f(a)\fL(S)/4}{2\fL(S)+2\fL(S)}\geq C_{11}.\ee
Note that (\ref{2.32}) also holds for $S\in \hB_3$. Using (\ref{3.25'}), (\ref{2.32}) and (\ref{2.48}),
we obtain
\be\label{2.48''}\ba{rl}&\P\(\Dp\bigcup_{S\in\hB_3}\{\fL(S,S^c)\geq {C\fL(S)}/4\ {\rm and\ }\fL(S)\geq |S|\}^c\)\\[5mm]
&\leq \Dp\sum_{s=an^d}^{(1-a)n^d}\b{n^d}{s}\exp(-\Omega (\sigma s\ln^{\zeta}n))\\[5mm]
&\leq \Dp\sum_{s=an^d}^{(1-a)n^d}\exp\(-\left[\Omega\(\sigma\ln^{\zeta}n\)-\ln(1/a)+1\right]s\)\rightarrow 0,\ {\rm as} \ n\rightarrow\infty.\ea\ee
 The lemma now follows immediately from (\ref{2.48'}) and (\ref{2.48''}).\QED

Now we can obtain the following lower bound for the edge isoperimetric constant $\iota$ of $G_n$.

\bp\label{p2.10} Suppose $0<\a<\beta<{1}/{2}$ with $\Gamma:=(2\beta)^d-(2\a)^d>{1}/{2}$. Then
\begin{description}
\item[] (i) for $0<\zeta<1$ and $\sigma>0$, or for $\zeta=0$ and $\sigma>0$ large enough, there exists $C_{12}>0$ such that \be\lim_{n\rightarrow\infty}\P\(\iota\geq {C_{12}}\ln^{-1}n\)=1;\ee
\item[] (ii) for $ \zeta>1 $ and $\sigma>0$, or for $\zeta=1$ and $\sigma>0$ large enough, there exists $C_{13}>0$ such that \be\lim_{n\rightarrow\infty}\P\(\iota\geq {C_{13}}\)=1.\ee
\end{description}
\ep

{\it Proof.} Let the following $\hB'_1$, $\hB'_2$ and $\hB'_3$ take the places of $\hB_1$, $\hB_2$ and $\hB_3$, the proposition follows from the same arguments as used in the proof of Proposition~\ref{p2.4}.
$$\ba{rl}&\hB'_1:=\{S\subset V_n: 1\leq |S|\leq M\ln n\},\ {\rm and}\\[3mm]&\hB'_2:=\{S\subset V_n: M\ln n<|S|\leq a|V_n|\},\ {\rm and}\\[3mm]
&\hB'_3:=\{S\subset V_n:a|V_n|<|S|\leq |V_n|/2\},\ea$$ where $M>0$ is a arbitrarily given constant and $a>0$ is given in Lemma~\ref{le2.3}.
\QED
\section{Proofs of Theorems~\ref{th1} and \ref{th2}}
\renewcommand{\theequation}{3.\arabic{equation}}
\setcounter{equation}{0}

{\it Proof of  A  of Theorem~\ref{th1}.} Noting that ${\rm diam}(G_n)$ is non-increasing in $\zeta$, we obtain part A of Theorem~\ref{th1} by Proposition~\ref{le2.2},  Proposition~\ref{p2.10} and the following Lemma~\ref{le3.1}.\QED

\bl\label{le3.1} For any connected graph $G=(V,E)$, let $\Delta(G)$ denote its maximum degree, $\iota(G)$ denote its edge isoperimetric constant and let ${\rm diam(G)}$ denote its diameter. Then \be\label{3.1}{\rm diam(G)}\leq \frac{4\Delta(G)}{\iota(G)}\ln |V|.\ee\el
{\it Proof.} This is a well known result in algebraic graph theory, for a detailed proof, one may refer to \cite{AM,BZ}.\QED

{\it Proofs of B of Theorem~\ref{th1} and Theorem~\ref{th2}.} Let $B_n$ is a box in $T^d_n$ with side length $2\ln^{r} n$, $r>0$. Let $A_n$ be the event that ${\cal L}(B_n,B_n^c)=0$. Then, for some $c_1>0$,
$$\P(A_n)=(1-p_n)^{(2\ln^{r} n)^{d}\Gamma n^d}\geq \exp\(-c_1\sigma\Gamma 2^d \ln^{dr+\zeta}n\),\ \rm{for\ large\ }n,$$
and
$$\ba{ll}&\P(\rm{there\  exists\ no\ box}\ B_n\ \rm{such\ that}\ A_n\ \rm{occurs})\\[2mm]
&\leq \left\{1-\exp\(-c_1\sigma\Gamma 2^d \ln^{dr+\zeta}n\)\right\}^{\(\frac{n}{2\ln^r n}\)^d}.\ea$$
By inequality $\ln(1-x)<-x$ for $0<x<1$, $\ln \P(\rm{there\  exists\ no\ box}\ B_n\ \rm{such\ that}\ A_n\ \rm{occurs})$ is bounded from above by
$$\Df{-n^d}{(2\ln^{r} n)^d}n^{-c_1\sigma 2^{d}\Gamma \ln^{\zeta+rd-1} n}.$$
Suppose $\zeta<1$, then for any $0<r<(1-\zeta)/d$, or for $r=(1-\zeta)/d$ and $\sigma>0$ small enough, the above quantity tends to $-\infty$ as $n\rightarrow\infty$. Hence
\be\label{3.3'}\lim_{n\rightarrow\infty}\P(\rm{there\ exists\ such\ a \ box}\ B_n\ \rm{such\ that}\ A_n \rm \ {occurs})=1.\ee
Still denote by $B_n$ the box with ${\cal L}(B_n,B_n^c)=0$. Note that the existence of such a box first implies that $\rm{diam}(G_n)\geq \ln^r n$, this finishes the proof of part B of Theorem~\ref{th1}.

 To finish the proof of part B of Theorem \ref{th2}, let's consider the lazy random walk $\{X_t\}$ starts at $o$, the center of the box $B_n$. By central limit Theorem, it is direct to check that, with probability tends to 1 (as $n\rightarrow\infty$), $\{X_t\}$ can not escape from $B_n$ in time $\ln^{\nu}n$ for any $\nu=2r<2(1-\zeta)/d$.
Then, for any $t\leq \ln^{\nu}n$, $P^t(o,u)=0$ for all $u\in V_n\setminus B_n$, so
$$||P^t(o,\cdot)-\pi||_{TV}\geq \frac 12\sum_{u\notin B_n}\pi(u)\geq\frac 12-\frac{|E(B_n)|}{|E_n|}\geq \frac 12-\frac{2d(2\ln^{r}n)^d}{dn^d}>\frac 1 e,\ \rm{for\ large}\ n.$$
Where $E(B_n)$ is the set of edges in $T^d_n$ with at least one endpoint in $B_n$, and $E_n$ is the edge set of $G_n$.
By definition of mixing time, we have $T_{\rm mix}\geq \ln^{\nu}n$
and finish the proof of B of Theorem~\ref{th2}.\QED

{\it Proof of A of Theorem~\ref{th2}.} For our lazy random walk $\{X_t:t\geq 0\}$ on $G_n$, matrix theory tell us that the transition kernel (P(u,v)) has nonnegative real eigenvalues $$1=\lambda_0\geq\lambda_1\geq\lambda_2\geq\ldots\geq\lambda_{n^d-1}\geq 0.$$ Note that $1-\lambda_1$ is called the spectral gap of $(P(u,v))$. Let $\pi_{\rm min}=\min_{u\in V_n}\pi(u)$.

 As a standard relation, it can be found in \cite[Theorem 12.5]{LPW} that
 \be\label{3.5}T_{\rm mix}\leq\ln\(\frac e{\pi_{\rm min}}\)\frac 1{1-\lambda_1}.\ee

The spectral gap $1-\lambda_1$ can be bounded from above and below by the conductance $h$ in the following way (see \cite[Theorem 6.2.1]{D}),
\be\label{3.6}\frac{h^2}2\leq 1-\lambda_1\leq 2h.\ee
On the other hand, it follows from Proposition~\ref{le2.2} that
\be\label{3.7}\P\(\Df 1{\pi_{\rm min}}\leq n^d\ln^{\zeta\vee 1} n\)\geq \P\(\Delta(G_n)\leq D\ln^{\zeta\vee 1} n\)\rightarrow 1, \ {\rm as\ }n\rightarrow\infty.\ee

Thus, the desired result follows from (\ref{2.22}), (\ref{3.5})-(\ref{3.7}) and Proposition~\ref{p2.4}.

\QED

\end{document}